\documentclass[english]{article}
\usepackage[T1]{fontenc}
\usepackage[latin9]{inputenc}
\usepackage{amsmath}
\usepackage{amssymb}
\usepackage{esint}
\usepackage{babel}
\usepackage{authblk}
\begin{document}

\title{Balanced Phase Field model for Active Surfaces}

\author[$1$]{Jozsef Molnar}
\author[$2$]{Peter Horvath}

\affil[$1$]{Synthetic and Systems Biology Unit\\
 Biological Research Centre,\\
 Hungarian Academy of Sciences\\
 Szeged, Hungary\\
 Email: jmolnar64@digikabel.hu}

\affil[$2$]{Synthetic and Systems Biology Unit\\
 Biological Research Centre,\\
 Hungarian Academy of Sciences\\
 Szeged, Hungary\\
 Email: horvath.peter@brc.mta.hu}

\maketitle
\begin{abstract}
In this paper we present a balanced phase field model for active surfaces.
This work is devoted to the generalization of the Balanced Phase Field
Model for Active Contours devised to eliminate the often undesirable
curvature-dependent shrinking of the zero level set while maintaining
the smooth interface necessary to calculate the fundamental geometric
quantities of the represented contour. As its antecedent work, the
proposed model extends the Ginzburg-Landau phase field energy with
a higher order smoothness term. The relative weights are determined
with the analysis of the level set motion in a curvilinear system
adapted to the zero level set. The proposed model exhibits strong
shape maintaining capability without signicant interference with the
active (e.g. a segmentation) model.
\end{abstract}

\section{Introduction}

Geometric active contours and surfaces \cite{caselles1993geometric,malladi1995shape}
are widely used for image segmentation where the representation of
contours/surfaces are mainly implicit: the zero level set of an appropriately
constructed function discretized on a fixed grid (Eulerian description).
The evolution of the level set function is governed by the Euler-Lagrange
equation associated with the appropriately designed functional for
the segmentation problem. Strict criteria are to be fulfilled by an
adequate level set representation. The most important one is that
it needs to be reasonably smooth across a certain neighborhood of
the zero level set to provide the basis of the accurate calculation
of fundamental geometric quantities of the contour/surface, the building
blocks of the equation(s) associated with the segmentation problem.
On the other hand, the segmentation equation deteriorates the shape
of the level set function - measurements must be taken to correct
it periodically. 

During the decades several methods were elaborated to cope with this
problem. The two main approaches are a) reinitialization and b) extension
of the PDE associated with the original problem with an extra term
that penalizes the deviations from the smooth (usually distance) function.
Reinitializing the level set function by calculating the distance
to the contours/surfaces on the whole domain is slow and may cause
instability at discontinuous locations of the distance function. The
partial remedy for this problems is the narrow band technique \cite{peng1999pde}
for the price of higher complexity. The extension of the original
PDE with a distance regularizing term \cite{li2010distance} may add
instability too (see \cite{zhang2013reinitialization}) or increase
complexity \cite{wang2014enhanced}\cite{wu2016indirectly}. More
importantly, these approaches may move the zero level set away from
the expected stopping location, which is rarely acceptable. From theoretical
perspective, any method dedicated to this shape maintaining should
have the least possible interference with the segmentation PDE.

The Ginzburg-Landau phase field model was introduced in the image
segmentation literature in \cite{rochery2005phase}\cite{horvath2007gas}.
It possesses interesting advantages over the earlier level set frameworks
as greater topological freedom; the possibility of a `neutral' initialization;
and a purely energy-based formulation. It also automatically forms
a narrow band around the zero level set with fast shape recovery owing
to a double well potential term incorporated to its functional; but
it still moves the level sets due to the energy proportional to the
length of the contour (or the surface area of the zero level set surface).
This problem was treated with high efficiency in {[}XXX{]} for active
contours. This work is devoted to the generalisation of the proposed
balanced phase field model for active surfaces.

The structure of the paper is the following. In section \ref{sec:Phase-field-models}
we summarize the Ginzburg-Landau and the balanced phase field model.
Then we examine the balanced phase field model for active surfaces
in section \ref{sec:The-balanced-phase-model}. Section \ref{sec:Discussion}
concludes the paper by discussion.

\section{Phase field models\label{sec:Phase-field-models}}

In the level set framework, the contours (2D) and surfaces (3D) are
represented by a constant (usually the zero) level set of a function
of two $\phi\left(x,y\right)$ and three variables $\Phi\left(x,y,z\right)$
respectively. The quantities of the segmentation problem are extracted
from these functions, such as the unit normal vector $\mathbf{n}=\frac{\nabla\phi}{\left|\nabla\phi\right|}$
or the curvature $\kappa=-\nabla\cdot\left(\frac{\nabla\phi}{\left|\nabla\phi\right|}\right)$
for contours and $\mathbf{n}=\frac{\nabla\Phi}{\left|\nabla\Phi\right|}$
or the sum curvarure $K_{S}=-\nabla\cdot\left(\frac{\nabla\Phi}{\left|\nabla\Phi\right|}\right)$
for surfaces, where $\nabla$ is the gradient operator of the appropriate
dimensions, ``$\cdot$'' stands for the scalar (dot) product, \textit{i.e.}
$\nabla\cdot\mathbf{v}$ is the divergence of the vector field $\mathbf{v}$.
The level set function is usually maintained on a uniform grid and
its derivatives are approximated by finite differences. This manner
of calculation requires the level set function to be approximately
linear locally, across a small neighborhood of the zero level set.
Phase field is one of the possible realizations of the level set frameworks.
Its energy functional is designed to form regions with $\pm1$ field
values (with the help of a double well potential term) and a smooth
transition between these regions adding smoothness term(s), naturally
representing a narrow band around the zero level set.

\subsection{Summary of the balanced phase field model for active contours}

The two dimensional balanced phase field model was introduced in {[}XXX{]}
with the aim to eliminate the undesired shrinking effect of the Ginzburg-Landau
phase field model. The proposed model extended the Ginzburg-Landau
phase field energy 
\begin{equation}
\iint_{\varOmega}\frac{D_{o}}{2}\left|\nabla\phi\right|^{2}+\lambda\left(\frac{\phi^{4}}{4}-\frac{\phi^{2}}{2}\right)dA\,.
\end{equation}
with a higher order smoothness term $\frac{D}{2}\left|\triangle\phi\right|^{2}$
(plus a constant term $\frac{\lambda}{4}$ such that the extended
functional expresses the energy of the transitional regions). The
relative weights were determined by the analysis of the extended energy
and the constant level set motion in a curvilinear system adapted
to the zero level set. Two conditions could be set: a) one for the
width of the transition (hereinafter denoted by $W$) between the
field values $\pm1$ (where both the Ginzburg-Landau and the balanced
functionals take their energy minima) and b) another for the elimination
of the curvature dependent term of the motion equation (i.e. the Euler-Lagrange
equation expressed in the adapted system) of the zero level set invoking
an adequately chosen ansatz. The approach led to two equations for
the weights in the extended functional as the functions of the width
of the transition. Note that the (minimal value of the) width required
is a priori known by the highest order of derivative occur in the
segmentation model.

Here we asses the most important results. Since two constraints have
to be satisfied, one of the weights can be arbitrarily set ($D_{o}$
is chosen to be $-1$). The balanced phase field functional and the
Euler-Lagrange equation then become
\begin{equation}
\iint_{\varOmega}\frac{W^{2}}{16}\left|\triangle\phi\right|^{2}-\frac{1}{2}\left|\nabla\phi\right|^{2}+\frac{21}{W^{2}}\left(\frac{\phi^{4}}{4}-\frac{\phi^{2}}{2}\right)dA
\end{equation}
and
\begin{equation}
\frac{W^{2}}{16}\triangle\triangle\phi+\triangle\phi+\frac{21}{W^{2}}\left(\phi^{3}-\phi\right)=0
\end{equation}
respectively. The gradient descent of the Euler-Lagrange equation
is recommended to be used for reinitialization using fix iteration
number $n\geq10$; with this value the balanced phase field is stable,
ensure smooth transition without significantly affecting the motion
of the constant level sets.

The question arises naturally: how can these results be extended for
the three dimensional case (for active surfaces).

\subsection{The three dimensional Ginzburg-Landau functional}

The energy of the simplest three dimensional Ginzburg-Landau phase
field level set representation: $\Phi\left(x,y,z\right)$ is defined
by the functional:
\begin{equation}
E\left(\Phi\right)=\iiint_{\Omega}\frac{1}{2}\left|\nabla\Phi\right|^{2}+\lambda\left(\frac{\Phi^{4}}{4}-\frac{\Phi^{2}}{2}+\frac{1}{4}\right)dxdydz,\label{eq:Energy_Cartesian-GiLa}
\end{equation}
where $\nabla\Phi$ is the gradient of the field $\Phi$, $\Omega$
represents the volume (the whole voxel image) of the integration.
The origin of the energy scale can be chosen freely. Term $\frac{\lambda}{4}$
is added such that at field values $\Phi=\pm1$ (where the functional
has its minima) $E=0$. This constant term does not influence the
Euler-Lagrange equation associated with the functional, which is:
\begin{equation}
-\triangle\Phi+\lambda\left(\Phi^{3}-\Phi\right)=0,\label{eq:Euler-Lagrange-GiLa}
\end{equation}
where $\triangle\Phi$ is the Laplacian of $\Phi$. As in the 2D case
it is easy to prove that energy (\ref{eq:Energy_Cartesian-GiLa})
is proportional to the surface area of the enclosed volume and as
a consequence the gradient descent of (\ref{eq:Euler-Lagrange-GiLa})
is driven by the sum curvature of the zero level set surface at every
point.

\section{The balanced phase field model for active surfaces\label{sec:The-balanced-phase-model}}

\subsection{The balanced phase field functional}

By analogy to the 2D version we propose the three dimensional balanced
phase field $\Phi\left(x,y,z\right)$ for level set representation
with the energy functional defined as:
\begin{equation}
E\left(\Phi\right)=\iiint_{\Omega}\frac{D}{2}\left|\triangle\Phi\right|^{2}-\frac{1}{2}\left|\nabla\Phi\right|^{2}+\lambda\left(\frac{\Phi^{4}}{4}-\frac{\Phi^{2}}{2}+\frac{1}{4}\right)dxdydz\,.\label{eq:Energy_Cartesian}
\end{equation}
Again, term $\frac{\lambda}{4}$ is added such that at field values
$\Phi=\pm1$ (where this functional has its minima) the energy becomes
zero and any deviation from the zero value is identified as the energy
of the transitional stripes between field values $-1$ and $1$. The
associated Euler-Lagrange equation is:
\begin{equation}
D\triangle\triangle\Phi+\triangle\Phi+\lambda\left(\Phi^{3}-\Phi\right)=0\,.\label{eq:Euler-Lagrange-1}
\end{equation}
We wish to determine the weights $D$ and $\lambda$ such that the
motion of the level sets governed by (\ref{eq:Euler-Lagrange-1})
are to be independent of the curvatures of the surfaces determined
by the level sets.

\subsection{The metric of the adapted system}

To get quantitative insight, we examine the system energy and the
motion of the zero level set in the curvilinear system adapted to
the zero level set. Let $\mathbf{S}\left(u,v\right)$ be the zero
level set surface, using Gaussian description. The space in the vicinity
of $\mathbf{S}$ can be parameterized as $\mathbf{R}\left(u,v,w\right)=\mathbf{S}\left(u,v\right)+w\mathbf{n}\left(u,v\right)$,
where $\mathbf{n}=\frac{\mathbf{S}_{u}\times\mathbf{S}_{v}}{\left|\mathbf{S}_{u}\times\mathbf{S}_{v}\right|}$
is the unit normal vector of the surface at point identified with
general coordinates $u,\, v$; lower indices stand for the partial
derivatives, \textit{i.e.} $\mathbf{S}_{u},\,\mathbf{S}_{v}$ are
the local (covariant) basis vectors. The length of the zero level
set surface normal vector $\left|\mathbf{S}_{u}\times\mathbf{S}_{v}\right|$
is equivalent to the square root of the determinant of the metric
tensor%
\footnote{Also known as first fundamental form.%
}: $\left[G_{ik}\right]=\left[\mathbf{S}_{i}\cdot\mathbf{S}_{k}\right]$
which is denoted by $\sqrt{G}$ (\textit{i.e.} $G=\det\left[G_{ik}\right]$).
It is used to define the parameterization independent infinitesimal
surface element $dS=\sqrt{G}dudv$. The square root of the determinant
of the metric tensor $\left[g_{ik}\right]=\left[\mathbf{R}_{i}\cdot\mathbf{R}_{k}\right]$,
$i,k\in\left\{ u,v,w\right\} $ is denoted by $\sqrt{g}$ and used
to define the parameterization independent infinitesimal volume element
$dV=\sqrt{g}dudvdw$. It can be expressed as the determinant of the
matrix constructed from the covariant basis vectors
\begin{eqnarray}
\mathbf{R}_{u} & = & \mathbf{S}_{u}+w\mathbf{n}_{u}\nonumber \\
\mathbf{R}_{v} & = & \mathbf{S}_{v}+w\mathbf{n}_{v}\\
\mathbf{R}_{w} & = & \mathbf{n}\nonumber 
\end{eqnarray}
such that $\sqrt{g}=\mathbf{n}\cdot\left(\mathbf{R}_{u}\times\mathbf{R}_{v}\right)$.
Expanding this expression we have:
\begin{alignat}{1}
\sqrt{g} & =\sqrt{G}\nonumber \\
 & +w\left[\mathbf{n}_{u}\cdot\left(\mathbf{S}_{v}\times\mathbf{n}\right)+\mathbf{n}_{v}\cdot\left(\mathbf{n}\times\mathbf{S}_{u}\right)\right]\label{eq:space-metric}\\
 & +w^{2}\left|\mathbf{n}_{u}\times\mathbf{n}_{v}\right|\,.\nonumber 
\end{alignat}
In the second line $\mathbf{S}_{v}\times\mathbf{n}=\sqrt{G}\mathbf{S}^{u}$,
$\mathbf{n}\times\mathbf{S}_{u}=\sqrt{G}\mathbf{S}^{v}$, where $\mathbf{S}^{u},$
$\mathbf{S}^{v}$ are the contravariant basis vectors of the surface
$\mathbf{S}$ with definition $\mathbf{S}^{i}\cdot\mathbf{S}_{k}=\delta_{k}^{i}$,
$i,k\in\left\{ u,v,w\right\} $ ($\delta_{k}^{i}$is the Cronecker
delta). The second line of (\ref{eq:space-metric}) is therefore the
$w\sqrt{G}$ times the divergence of the unit normal vector which
is in turn the negative of the sum cirvature $-K_{S}$. $\left|\mathbf{n}_{u}\times\mathbf{n}_{v}\right|$
in the third line is the integrand of the total curvature expression
equivalent with $\sqrt{G}K_{G}$ where $K_{G}$ is the Gaussian curvature
(see also appendix A). The square root of the metric therefore can
be expressed by a quadratic function of $w$ with coefficients being
the sum and Gaussian curvatures of the zero level set:
\begin{equation}
\sqrt{g}=\sqrt{G}\left(1-wK_{S}+w^{2}K_{G}\right)\,.\label{eq:metric}
\end{equation}

\subsection{Energy terms in the adapted system\label{sub:The-energy-constituents}}

First we examine the constituents of energy (\ref{eq:Energy_Cartesian})
in the curvilinear system adatpted to the level sets surfaces. In
this case $\Phi\left(u,v,w\right)$ takes constant values regardless
the parameter values $u,\, v$, hence its partial derivatives wrt
these parameters are all zero, that is $\frac{\partial^{n}\Phi}{\partial u^{n}}=0$,
$\frac{\partial^{m}\Phi}{\partial v^{m}}=0$, $m,\, n$ are arbitrary.
The gradient 
\begin{equation}
\left.\frac{\partial\Phi}{\partial u}\mathbf{R}^{u}+\frac{\partial\Phi}{\partial v}\mathbf{R}^{v}+\frac{\partial\Phi}{\partial w}\mathbf{n}\right|_{\Phi=const}=\frac{\partial\Phi}{\partial w}\mathbf{n},
\end{equation}
and the Laplacian
\begin{eqnarray}
\left.\frac{1}{\sqrt{g}}\frac{\partial\left(\sqrt{g}g^{ik}\frac{\partial\Phi}{\partial u^{k}}\right)}{\frac{\partial\Phi}{\partial u^{i}}}\right|_{\Phi=const} & = & \frac{\partial^{2}\Phi}{\partial w^{2}}+\frac{1}{\sqrt{g}}\frac{\partial\sqrt{g}}{\partial w}\frac{\partial\Phi}{\partial w},\label{eq:Laplacian_adapted}
\end{eqnarray}
where
\begin{equation}
\frac{1}{\sqrt{g}}\frac{\partial\sqrt{g}}{\partial w}=\frac{-K_{S}+wK_{G}}{1-wK_{S}+w^{2}K_{G}}\,.\label{eq:gene_curve}
\end{equation}
Note that in the general expression (left of (\ref{eq:Laplacian_adapted}))
the Einstein summation convention is used. These expressions are dependent
only on the geometric quantities of the zero level set $K_{S},\, K_{G}$
and the derivatives of the level set function in normal direction.
To simplify the notation, from now on we use primes to denote the
derivatives in the normal direction: $\Phi^{\prime}\equiv\frac{\partial\Phi}{\partial w}$,
$\Phi^{\prime\prime}\equiv\frac{\partial^{2}\Phi}{\partial w^{2}}$...;
notice that both the gradient and the Laplacian expressions contains
derivatives only in normal direction explicitely. Implicitely the
derivatives wrt $u$ and $v$ occur in the geometric quantities of
the zero level set surface only ($K_{S}$, $K_{G}$).

At this point it is tempting to assume the following 
\begin{enumerate}
\item $\Phi\left(u,v,w\right)\equiv\Phi\left(w\right)$, that is the constant
level sets are equidistant to each-other (hence $\frac{\partial^{l+m+n}\Phi}{\partial u^{l}\partial v^{m}\partial w^{n}}\equiv0$); 
\item They are arranged symmetrically around the zero level set $\Phi\left(0\right)$
i.e. $\Phi\left(-w\right)=-\Phi\left(w\right)$; 
\item The transition between $\Phi=-1$, $\Phi=1$ (representing the minimal
energy states) is confined to a stripe with constant width $W$. 
\end{enumerate}
Note that from condition 2. $\Phi\left(0\right)=0$. These assumptions
are certainly true for the plane (for symmetry reason) and violated
only wherever curvatures are present; for this reason the low curvature
condition 
\begin{equation}
1-wK_{S}+w^{2}K_{G}\approx1\label{eq:curvexpress_approx}
\end{equation}
needs to be assumed.

\subsection{Ansatz for the level set function}

The simplest possible ansatz satisfying the assumption taken in \ref{sub:The-energy-constituents}
is the cubic function $\Phi\doteq aw^{3}+bw$ with boundary conditions:
\begin{eqnarray}
\Phi\left(\frac{W}{2}\right) & = & 1\nonumber \\
\Phi\left(-\frac{W}{2}\right) & = & -1\label{eq:ansatz_prop}\\
\Phi^{\prime}\left(\pm\frac{W}{2}\right) & = & 0\,.\nonumber 
\end{eqnarray}
With these, the function and its derivatives involved in the system
energy are:
\begin{eqnarray}
\Phi\left(w\right) & =- & \frac{4}{W^{3}}w^{3}+\frac{3}{W}w\nonumber \\
\Phi^{\prime}\left(w\right) & =- & \frac{12}{W^{3}}w^{2}+\frac{3}{W}\label{eq:ansatz}\\
\Phi^{\prime\prime}\left(w\right) & =- & \frac{24}{W^{3}}w\nonumber 
\end{eqnarray}

\subsection{Energy expression in the adapted system}

With the one-dimensional ansatz, the energy (\ref{eq:Energy_Cartesian})
becomes:

\begin{equation}
\varoiint\int\left[\frac{D}{2}\left(\Phi^{\prime\prime}+\frac{1}{\sqrt{g}}\frac{\partial\sqrt{g}}{\partial w}\Phi^{\prime}\right)^{2}-\frac{1}{2}\left(\Phi^{\prime}\right)^{2}+\lambda\left(\frac{\Phi^{4}}{4}-\frac{\Phi^{2}}{2}+\frac{1}{4}\right)\right]\sqrt{g}dwdudv\,.\label{eq:Energy_Adapted}
\end{equation}
Substituting the ansatz (\ref{eq:ansatz}) into energy (\ref{eq:Energy_Adapted})
and using the low-curvature approximation (\ref{eq:curvexpress_approx}),
the energy, as the function of the width of the transition, is
\begin{equation}
\tilde{E}\left(W\right)=\left(D\frac{24}{W^{3}}-\frac{12}{5W}+\lambda\frac{W}{10}\right)A\label{eq:energy_by_width}
\end{equation}
(see appendices B, C and D).

\subsection{Optimal width}

Handling the energy expression (\ref{eq:energy_by_width}) as extreme
value problem%
\footnote{This is rational, because the zero level set surface practically static
(does not move) during the time necessary to form the shape of transition.%
} one can get an equation for the width of the transition as the function
of two parameters - the weights of the constituents of (\ref{eq:Euler-Lagrange-1}):
\begin{equation}
\frac{d\tilde{E}}{dW}=\left(-\frac{72D}{W^{4}}+\frac{12}{5W^{2}}+\frac{\lambda}{10}\right)A\doteq0\,.
\end{equation}
Rearranging wherever surface area is not zero ($A>0$) we obtain to
the first equation we need:
\begin{multline}
\lambda W^{4}+24W^{2}-720D=0\:\rightarrow\\
\qquad W^{2}=\frac{1}{\lambda}\left(-12\pm\sqrt{12^{2}+720D\lambda}\right)\,.\label{eq:first-equation}
\end{multline}

\subsection{Euler-Lagrange equation in the adapted system}

Under the three conditions stated in point \ref{sub:The-energy-constituents},
the (approximate) Euler-Lagrange equation (\ref{eq:Euler-Lagrange-1})
in the adapted system becomes a fourth order ordinary differential
equation:
\begin{multline}
D\left\{ \frac{\partial^{4}\Phi}{\partial w^{4}}+\frac{2}{\sqrt{g}}\frac{\partial\sqrt{g}}{\partial w}\frac{\partial^{3}\Phi}{\partial w^{3}}+\left[\left(\frac{1}{\sqrt{g}}\frac{\partial\sqrt{g}}{\partial w}\right)^{2}+2\frac{\partial}{\partial w}\left(\frac{1}{\sqrt{g}}\frac{\partial\sqrt{g}}{\partial w}\right)\right]\frac{\partial^{2}\Phi}{\partial w^{2}}\right.\\
\left.\qquad\qquad+\left[\triangle_{T}\left(\frac{1}{\sqrt{g}}\frac{\partial\sqrt{g}}{\partial w}\right)+\left(\frac{\partial^{2}}{\partial w^{2}}+\frac{1}{\sqrt{g}}\frac{\partial\sqrt{g}}{\partial w}\frac{\partial}{\partial w}\right)\left(\frac{1}{\sqrt{g}}\frac{\partial\sqrt{g}}{\partial w}\right)\right]\frac{\partial\Phi}{\partial w}\right\} \\
+\frac{\partial^{2}\Phi}{\partial w^{2}}+\frac{1}{\sqrt{g}}\frac{\partial\sqrt{g}}{\partial w}\frac{\partial\Phi}{\partial w}+\lambda\left(\Phi^{3}-\Phi\right)=0\label{eq:Euler-Lagrange-2}
\end{multline}
(see appendix E).

\subsection{Motion of the level sets}

According to the equidistance condition - assumed to be persistent
during the evolution governed by the Euler-Lagrange equation (\ref{eq:Euler-Lagrange-1})
(or approximate equation (\ref{eq:Euler-Lagrange-2})), it is sufficient
to examine the motion of any constant level set. The simplest case
is the zero level set; wrt this set the antisymmetry condition $\Phi\left(w\right)-\Phi\left(-w\right)=0$
(assumed in point \ref{sub:The-energy-constituents}) and consequently
$\frac{\partial^{2k}\Phi}{\partial w^{2k}}$, $k=1,2...$ are satisfied.
Moreover, from (\ref{eq:gene_curve}):
\begin{eqnarray}
\left.\frac{1}{\sqrt{g}}\frac{\partial\sqrt{g}}{\partial w}\right|_{w=0} & = & -K_{S}\nonumber \\
\left.\left(\frac{\partial^{2}}{\partial w^{2}}+\frac{1}{\sqrt{g}}\frac{\partial\sqrt{g}}{\partial w}\frac{\partial}{\partial w}\right)\left(\frac{1}{\sqrt{g}}\frac{\partial\sqrt{g}}{\partial w}\right)\right|_{w=0} & = & -K_{S}\left(K_{S}^{2}-4K_{G}\right)\label{eq:zero-set-quant}\\
\left.\triangle_{T}\left(\frac{1}{\sqrt{g}}\frac{\partial\sqrt{g}}{\partial w}\right)\right|_{w=0} & = & -\triangle_{T}K_{S}\nonumber 
\end{eqnarray}
where $\triangle_{T}$ denotes the tangential components of the Laplace
operator. Substituting (\ref{eq:zero-set-quant}) to (\ref{eq:Euler-Lagrange-2}),
the adapted Euler-Lagrange equation for the zero level set is reduced
to:
\begin{multline}
-2DK_{S}\frac{\partial^{3}\Phi}{\partial w^{3}}-D\left[\triangle_{T}K_{S}+K_{S}\left(K_{S}^{2}-4K_{G}\right)\right]\frac{\partial\Phi}{\partial w}-K_{S}\frac{\partial\Phi}{\partial w}=0\,.\label{eq:Euler-Lagrange-zeroset}
\end{multline}
The sum curvature dependency is therefore can be eliminated by the
condition (involving the 1st and the 3rd term in (\ref{eq:Euler-Lagrange-zeroset})):
\begin{equation}
\left.-2\frac{\partial^{3}\Phi}{\partial w^{3}}-\frac{\partial\Phi}{\partial w}\right|_{w=0}\doteq0,
\end{equation}
or using ansatz (\ref{eq:ansatz}):
\begin{equation}
2D\frac{24}{W^{3}}-\frac{3}{W}\doteq0\,.\label{eq:second_equation}
\end{equation}

\subsection{Energy expression and Euler-Lagrange equation for curvature-independent
motion}

Equations (\ref{eq:first-equation}) and (\ref{eq:second_equation})
determines parameters $\lambda$ and $D$ in energy (\ref{eq:Energy_Cartesian})
with the curvature-driven shrinking effect removed from the gradient
descent of its associated Euler-Lagrange equation (\ref{eq:Euler-Lagrange-1})
as the function of the width of transition $W$. The solution is:
\begin{eqnarray}
D\left(W\right) & = & \frac{W^{2}}{16}\nonumber \\
\lambda\left(W\right) & = & \frac{21}{W^{2}}\,.\label{eq:weights}
\end{eqnarray}
Note that the gradient descent of the balanced zero level set equation
- the 2nd term in (\ref{eq:Euler-Lagrange-zeroset}) - still describe
dynamic surface, but with a motion of a very modest pace. In fact
the remaining term $\triangle_{T}K_{S}+K_{S}\left(K_{S}^{2}-4K_{G}\right)$
is very close to the solution $\triangle_{T}K_{S}+\frac{1}{2}K_{S}\left(K_{S}^{2}-4K_{G}\right)$
of the functional derivative associated with the Euler's elastica
of surfaces $\frac{1}{2}\oiint K_{S}^{2}dA$.

With the determined weights, we have the energy (\ref{eq:Energy_Cartesian})
\begin{equation}
E\left(\Phi\right)=\iiint_{\Omega}\frac{W^{2}}{32}\left|\triangle\Phi\right|^{2}-\frac{1}{2}\left|\nabla\Phi\right|^{2}+\frac{21}{W^{2}}\left(\frac{\Phi^{4}}{4}-\frac{\Phi^{2}}{2}+\frac{1}{4}\right)dxdydz\label{eq:Balanced-energy-Cartesian}
\end{equation}
and the Euler-Lagrange equation associated with it
\begin{equation}
\frac{W^{2}}{16}\triangle\triangle\Phi+\triangle\Phi+\frac{21}{W^{2}}\left(\Phi^{3}-\Phi\right)=0\,.\label{eq:Balanced-E-L-Cartesian}
\end{equation}

\section{Discussion\label{sec:Discussion}}

In this paper we generalized the 2D balanced field model to active
surfaces. It is shown by the examination of the equations of motions
- the Euler-Lagrange equations associated with the Ginzburg-Landau
and the balanced phase field models in the adapted curvilinear systems
- that (as usual) the sum curvature for active surfaces has the same
role as the curvature for active contours and can be eliminated using
the same constraints. This curvature/sum curvature correspondence
holds for the constraints that can be imposed todetermine the optimal
widths of the transitions. 

We concluded that the 3D equations expressed in Cartesian coordinates
have exactly same form as their 2D counterparts. As in 2D, the gradient
descent of the proposed model exhibits very fast shape recovery without
moving the zero level set significantly. In fact the motion of level
sets is similar to the motion associated with the Euler's elastica.
This remaining term contains a nonlinear expression of the sum and
Gaussian curvatures (expressible with a cubic polynomial of the principal
curvatures) and under the low curvature assuption its interference
with the segmentation model is negligible, the property that makes
this level set formulation suitable for accurate segmentation. As
in 2D, this balancing could be used for any model that includes Laplacian
smoothness term in their gradient descent equation like the reaction-diffusion
model.

\part*{\newpage{}Appendices}

\subsection*{Appendix A: The Gaussian term of the metric\label{sub:A:-The-Gaussian}}

The invariant surface element is defined with the infinitesimal area
of parallelogram spanned by the covariant basis vectors $\mathbf{S}_{u},\,\mathbf{S}_{v}$
as $dS=\left|\mathbf{S}_{u}\times\mathbf{S}_{v}\right|dudv$, where
the factor $\left|\mathbf{S}_{u}\times\mathbf{S}_{v}\right|$ is:
\begin{eqnarray}
\left|\mathbf{S}_{u}\times\mathbf{S}_{v}\right| & = & \sqrt{\left(\mathbf{S}_{u}\times\mathbf{S}_{v}\right)\cdot\left(\mathbf{S}_{u}\times\mathbf{S}_{v}\right)}\nonumber \\
 & = & \sqrt{\mathbf{S}_{u}\cdot\left[\mathbf{S}_{v}\times\left(\mathbf{S}_{u}\times\mathbf{S}_{v}\right)\right]}\nonumber \\
 & = & \sqrt{\mathbf{S}_{u}\cdot\left[\mathbf{\mathbf{\left(\mathbf{S}_{v}\cdot\mathbf{S}_{v}\right)S}_{u}-\left(\mathbf{S}_{u}\cdot\mathbf{S}_{v}\right)S}_{v}\right]}\label{eq:metricequivalence}\\
 & = & \sqrt{\left(\mathbf{S}_{u}\cdot\mathbf{S}_{u}\right)\left(\mathbf{S}_{v}\cdot\mathbf{S}_{v}\right)-\left(\mathbf{S}_{u}\cdot\mathbf{S}_{v}\right)^{2}}=\sqrt{G}\,.\nonumber 
\end{eqnarray}
(For the derivation, the triple scalar product $\mathbf{a}\cdot\left(\mathbf{b}\times\mathbf{c}\right)=\mathbf{b}\cdot\left(\mathbf{c}\times\mathbf{a}\right)$
and the triple cross product $\mathbf{a}\times\left(\mathbf{b}\times\mathbf{c}\right)=\left(\mathbf{a}\cdot\mathbf{c}\right)\mathbf{b}-\left(\mathbf{a}\cdot\mathbf{b}\right)\mathbf{c}$
equivalences are used.)

The partial derivatives of the unit normal vector $\mathbf{n}_{u}$
and $\mathbf{n}_{v}$ are the elements of the tangent space hence
can be decomposed such that $\mathbf{n}_{i}=\left(\mathbf{n}_{i}\cdot\mathbf{S}^{u}\right)\mathbf{S}_{u}+\left(\mathbf{n}_{i}\cdot\mathbf{S}^{v}\right)\mathbf{S}_{v}$,
$i\in\left\{ u,v\right\} $, where $\mathbf{S}_{u},\,\mathbf{S}_{v}$
are the local (covariant) basis, $\mathbf{S}^{u},\,\mathbf{S}^{v}$
are the contravariant basis vectors with the property $\mathbf{S}^{i}\cdot\mathbf{S}_{k}=\delta_{k}^{i}$
($\delta_{k}^{i}$ is the Kronecker delta). It can be seen by simple
substitution that
\begin{eqnarray}
\mathbf{S}^{u} & = & \frac{1}{\left|\mathbf{S}_{u}\times\mathbf{S}_{v}\right|}\mathbf{S}_{v}\times\mathbf{n}\nonumber \\
\mathbf{S}^{v} & = & \frac{1}{\left|\mathbf{S}_{u}\times\mathbf{S}_{v}\right|}\mathbf{n}\times\mathbf{S}_{u}\,.\label{eq:contra-by-co}
\end{eqnarray}
We also need
\begin{eqnarray}
\mathbf{n}\cdot\mathbf{S}_{k}\equiv0 & \rightarrow & \mathbf{n}_{i}\cdot\mathbf{S}_{k}=-\mathbf{n}\cdot\mathbf{S}_{ik}\nonumber \\
 &  & i,k\in\left\{ u,v\right\} \label{eq:to-secfundform}
\end{eqnarray}
Next we calculate $\mathbf{n}_{u}\times\mathbf{n}_{v}$ as
\begin{alignat}{1}
\mathbf{n}_{u}\times\mathbf{n}_{v}= & \left[\left(\mathbf{n}_{u}\cdot\mathbf{S}^{u}\right)\mathbf{S}_{u}+\left(\mathbf{n}_{u}\cdot\mathbf{S}^{v}\right)\mathbf{S}_{v}\right]\times\left[\left(\mathbf{n}_{v}\cdot\mathbf{S}^{u}\right)\mathbf{S}_{u}+\left(\mathbf{n}_{v}\cdot\mathbf{S}^{v}\right)\mathbf{S}_{v}\right]\nonumber \\
= & \left(\mathbf{S}^{u}\times\mathbf{S}^{v}\right)\left[\left(\mathbf{n}_{u}\cdot\mathbf{S}_{u}\right)\left(\mathbf{n}_{v}\cdot\mathbf{S}_{v}\right)-\left(\mathbf{n}_{u}\cdot\mathbf{S}_{v}\right)^{2}\right]\\
= & \left(\mathbf{S}^{u}\times\mathbf{S}^{v}\right)\left[\left(\mathbf{n}\cdot\mathbf{S}_{uu}\right)\left(\mathbf{n}\cdot\mathbf{S}_{vv}\right)-\left(\mathbf{n}\cdot\mathbf{S}_{uv}\right)^{2}\right]\,.\nonumber 
\end{alignat}
In the third line (\ref{eq:to-secfundform}) is used. Substitution
of (\ref{eq:contra-by-co}) leads to
\begin{alignat}{1}
\mathbf{n}_{u}\times\mathbf{n}_{v}= & \left(\mathbf{S}_{v}\times\mathbf{n}\right)\left(\mathbf{n}\times\mathbf{S}_{u}\right)\frac{\left(\mathbf{n}\cdot\mathbf{S}_{uu}\right)\left(\mathbf{n}\cdot\mathbf{S}_{vv}\right)-\left(\mathbf{n}\cdot\mathbf{S}_{uv}\right)^{2}}{\left|\mathbf{S}_{u}\times\mathbf{S}_{v}\right|^{2}}\nonumber \\
= & \mathbf{n}\left[\mathbf{S}_{u}\cdot\left(\mathbf{S}_{v}\times\mathbf{n}\right)\right]\frac{\det\left[\mathbf{II}\right]}{\det\left[G_{ik}\right]}\label{eq:Gaussian-metric-vector}\\
= & \left|\mathbf{S}_{u}\times\mathbf{S}_{v}\right|K_{G}\mathbf{n},\nonumber 
\end{alignat}
where the Gaussian curvature is given as the ratio of the determinants
of the second and first fundamental forms. The last line of equation
(\ref{eq:Gaussian-metric-vector}) is a vector with length
\begin{equation}
\left|\mathbf{n}_{u}\times\mathbf{n}_{v}\right|=\left|\mathbf{S}_{u}\times\mathbf{S}_{v}\right|K_{G}=\sqrt{G}K_{G}\,.\label{eq:Gaussian-metric-term}
\end{equation}

\subsection*{Appendix B: Approximation of the gradient integral term\label{sub:B:-Gradient-term}}

The second term of (\ref{eq:Energy_Adapted}) is
\begin{gather}
\varoiint\int-\frac{1}{2}\left(\Phi^{\prime}\right)^{2}\left(1-wK_{S}+w^{2}K_{G}\right)dwdA,
\end{gather}
where $dA=\sqrt{G}dudv$. Substituting the corresponding ansatz (\ref{eq:ansatz})
$\left(\Phi^{\prime}\right)^{2}=\frac{12^{2}}{W^{6}}w^{4}-\frac{72}{W^{4}}w^{2}+\frac{9}{W^{2}}$
and performing the integration between the boundaries $-\frac{W}{2},\frac{W}{2}$,
the result is:
\begin{equation}
-\frac{12}{5W}A-\frac{W^{3}}{24}\varoiint\frac{1}{R_{1}}\frac{1}{R_{2}}dA\,.\label{eq:grad-term-integ}
\end{equation}
where $A=\varoiint dA$ is the surface area of the zero level set
(assumed to be closed, hence the notation $\oiint$). The second term
is: 
\begin{equation}
\frac{W^{3}}{24}\varoiint\frac{1}{R_{1}}\frac{1}{R_{2}}dA=\frac{W}{24}\varoiint\frac{W}{R_{1}}\frac{W}{R_{2}}dA\ll\frac{W}{24}A,
\end{equation}
where $R_{1},R_{2}$ are the rays of the osculating circles in the
principal directions. If the ratio of $\frac{WA}{24}$ and the first
$\frac{12A}{5W}$ terms ($\thickapprox0.132W$) in (\ref{eq:grad-term-integ})
is not extremely big (a case fur moderate curvature valus) then the
second term can be omitted from (\ref{eq:grad-term-integ}).

\subsection*{Appendix C: Approximation of the Laplacian integral term\label{sub:C:-Laplacian-term}}

Here we calculate the first term of (\ref{eq:Energy_Adapted}) using
the cubic ansatz (\ref{eq:ansatz}), the metric (\ref{eq:metric})
and the invariant surface element expression $dA=\sqrt{G}dudv$. 
\begin{gather}
\varoiint\int\left(\Phi^{\prime\prime}+\frac{-K_{S}+2wK_{G}}{1-wK_{S}+w^{2}K_{G}}\Phi^{\prime}\right)^{2}\left(1-wK_{S}+w^{2}K_{G}\right)dwdA=\nonumber \\
\varoiint\int\left(1-wK_{S}+w^{2}K_{G}\right)\left(\Phi^{\prime\prime}\right)^{2}\nonumber \\
+2\left(-K_{S}+2wK_{G}\right)\Phi^{\prime\prime}\Phi^{\prime}\\
+\frac{\left(-K_{S}+2wK_{G}\right)^{2}}{1-wK_{S}+w^{2}K_{G}}\left(\Phi^{\prime}\right)^{2}dwdA\nonumber 
\end{gather}
In the second term, $2\Phi^{\prime\prime}\Phi^{\prime}=\left(\left(\Phi^{\prime}\right)^{2}\right)^{\prime}$.
Applying integration by parts to this integrand, on of the term $\left[\left(-K_{S}+2wK_{G}\right)\Phi^{\prime}\right]_{-\frac{W}{2}}^{\frac{W}{2}}=0$
(according to the third property (\ref{eq:ansatz_prop}) of the chosen
ansatz). What remains is:
\begin{gather}
\varoiint\int\left(1-wK_{S}+w^{2}K_{G}\right)\left(\Phi^{\prime\prime}\right)^{2}\nonumber \\
+\left[\frac{\left(-K_{S}+2wK_{G}\right)^{2}}{1-wK_{S}+w^{2}K_{G}}-2K_{G}\right]\left(\Phi^{\prime}\right)^{2}dwdA=\nonumber \\
\varoiint\int\left(1-wK_{S}+w^{2}K_{G}\right)\left(\Phi^{\prime\prime}\right)^{2}\\
+\frac{K_{S}^{2}+2K_{G}\left(w^{2}K_{G}-wK_{S}-1\right)}{1-wK_{S}+w^{2}K_{G}}\left(\Phi^{\prime}\right)^{2}dwdA\,.\nonumber 
\end{gather}
Now we use approximation (\ref{eq:curvexpress_approx}) and arrive
to:
\begin{gather}
\varoiint\int\left(\Phi^{\prime\prime}\right)^{2}+\left(K_{S}^{2}-2K_{G}\right)\left(\Phi^{\prime}\right)^{2}dwdA=\nonumber \\
\varoiint\int\left(\Phi^{\prime\prime}\right)^{2}+\left(K_{1}^{2}+K_{2}^{2}\right)\left(\Phi^{\prime}\right)^{2}dwdA,
\end{gather}
where $K_{1}$ and $K_{2}$ are the principal curvatures. With the
ansatz values (\ref{eq:ansatz}) the first term of the energy (\ref{eq:Laplacian_adapted})
integrated between the boundary values $-\frac{W}{2},\frac{W}{2}$
becomes
\begin{equation}
D\left(\frac{24}{W^{3}}A+\frac{12}{5W}\oiint\frac{1}{R_{1}^{2}}+\frac{1}{R_{2}^{2}}dA\right)\,.\label{eq:Lapl-term-integ}
\end{equation}
The second term is: 
\begin{equation}
\frac{12}{5W}\oiint\frac{1}{R_{1}^{2}}+\frac{1}{R_{2}^{2}}dA=\frac{12}{5W^{3}}\oiint\left(\frac{W}{R_{1}}\right)^{2}+\left(\frac{W}{R_{2}}\right)^{2}dA\ll\frac{12}{5W^{3}}A\,.
\end{equation}
Observing that the ratio of $\frac{12A}{5W^{3}}$ and the first $\frac{24A}{W^{3}}$
terms ($=0.5$) is not an extremely big number, one can conclude that
the second term can be omitted from (\ref{eq:Lapl-term-integ}).

\subsection*{Appendix D: Approximation of the double potential well term\label{sub:D:-Potentialwell-term}}

The 
\begin{multline}
\varoiint\int\lambda\left(\frac{\Phi^{4}}{4}-\frac{\Phi^{2}}{2}+\frac{1}{4}\right)\sqrt{g}dwdudv=\\
\qquad\varoiint\int\lambda\left(\frac{\Phi^{4}}{4}-\frac{\Phi^{2}}{2}+\frac{1}{4}\right)\left(1-wK_{S}+w^{2}K_{G}\right)dwdA
\end{multline}
 can be considered as the active term of the balanced phase field
model. Substituting the ansatz (\ref{eq:ansatz}) and carrying out
the integration between the boundary values $-\frac{W}{2},\frac{W}{2}$
one can have:
\begin{equation}
\lambda\left(\frac{W}{10}A+\frac{W^{3}}{12}\varoiint\frac{1}{R_{1}}\frac{1}{R_{2}}dA\right)=\lambda W\left(\frac{1}{10}A+\frac{1}{12}\varoiint\frac{W}{R_{1}}\frac{W}{R_{2}}dA\right)\,.
\end{equation}
Here the approximation $\int\frac{\Phi^{4}}{4}-\frac{\Phi^{2}}{2}+\frac{1}{4}dw\thickapprox0.1W$
is used. With similar reasoning to the previous cases it is obvious
that the second term can be ignored similarly. Note that the terms
neglegted in comparision with the dominant terms in the (gradient,
Laplacian and active part) one by one. If there is magnitudes of differences
between the relative weights, then this approximations can be invalid.
The final result however justified this approach a posteriory for
a wide range of realistic width selection.

\subsection*{Appendix E: Approximation of Laplacian and the Laplacian of Laplacian\label{sub:E:-LaplLapl-term}}

It is expedient to decompose the Lagrangian operator in the direction
tangential and perpendicular to the level sets
\begin{equation}
\triangle=\frac{1}{\sqrt{g}}\frac{\partial\left(\sqrt{g}g^{ik}\frac{\partial}{\partial u^{k}}\right)}{\frac{\partial}{\partial u^{i}}}\doteq\triangle_{T}+\left(\frac{\partial^{2}}{\partial w^{2}}+\frac{1}{\sqrt{g}}\frac{\partial\sqrt{g}}{\partial w}\frac{\partial}{\partial w}\right),
\end{equation}
$i,k\in\left\{ u,v,w\right\} $, the metric $\sqrt{g}$ is given by
(\ref{eq:metric}) and
\begin{multline}
\triangle_{T}=g^{uu}\frac{\partial^{2}}{\partial u^{2}}+2g^{uv}\frac{\partial^{2}}{\partial u\partial v}+g^{vv}\frac{\partial^{2}}{\partial v^{2}}\\
\quad+\frac{1}{\sqrt{g}}\left(\frac{\partial g^{uu}\sqrt{g}}{\partial u}+\frac{\partial g^{uv}\sqrt{g}}{\partial v}\right)\frac{\partial}{\partial u}\\
+\frac{1}{\sqrt{g}}\left(\frac{\partial g^{uv}\sqrt{g}}{\partial u}+\frac{\partial g^{vv}\sqrt{g}}{\partial v}\right)\frac{\partial}{\partial v}\,.
\end{multline}
It immediately follows that on level sets $\Phi=const$, the derivatives
$\frac{\partial^{m+n}\Phi}{\partial u^{m}\partial v^{n}}$ are automatically
zero, the tangential component $\triangle_{T}$ has no effect, i.e.
\begin{equation}
\triangle\Phi=\frac{\partial^{2}\Phi}{\partial w^{2}}+\frac{1}{\sqrt{g}}\frac{\partial\sqrt{g}}{\partial w}\frac{\partial\Phi}{\partial w}\,.\label{eq:lapl}
\end{equation}

Applying the Laplace operator twice leads to the decomposition
\begin{multline}
\triangle\triangle\Phi=\triangle_{T}\triangle_{T}\Phi+\triangle_{T}\left(\frac{\partial^{2}\Phi}{\partial w^{2}}+\frac{1}{\sqrt{g}}\frac{\partial\sqrt{g}}{\partial w}\frac{\partial\Phi}{\partial w}\right)\\
+\left(\frac{\partial^{2}}{\partial w^{2}}+\frac{1}{\sqrt{g}}\frac{\partial\sqrt{g}}{\partial w}\frac{\partial}{\partial w}\right)\triangle_{T}\Phi\\
\qquad\qquad\qquad\qquad+\left(\frac{\partial^{2}}{\partial w^{2}}+\frac{1}{\sqrt{g}}\frac{\partial\sqrt{g}}{\partial w}\frac{\partial}{\partial w}\right)\left(\frac{\partial^{2}\Phi}{\partial w^{2}}+\frac{1}{\sqrt{g}}\frac{\partial\sqrt{g}}{\partial w}\frac{\partial\Phi}{\partial w}\right)\,.\label{eq:lapl-lapl}
\end{multline}
The first term is automatically zero. If the level sets are equdistant
to each-other (i.e. $\frac{\partial^{m+n+r}\Phi}{\partial u^{m}\partial v^{n}\partial w^{r}}=0$)
then the non-zero terms remained are
\begin{multline}
\triangle\triangle\Phi=\frac{\partial\Phi}{\partial w}\triangle_{T}\left(\frac{1}{\sqrt{g}}\frac{\partial\sqrt{g}}{\partial w}\right)\\
\qquad\qquad+\left(\frac{\partial^{2}}{\partial w^{2}}+\frac{1}{\sqrt{g}}\frac{\partial\sqrt{g}}{\partial w}\frac{\partial}{\partial w}\right)\left(\frac{\partial^{2}\Phi}{\partial w^{2}}+\frac{1}{\sqrt{g}}\frac{\partial\sqrt{g}}{\partial w}\frac{\partial\Phi}{\partial w}\right),\label{eq:lapl-lapl-1}
\end{multline}
or rearranging the equation by the orders of the derivatives of the
phase field function
\begin{multline}
\triangle\triangle\Phi=\frac{\partial^{4}\Phi}{\partial w^{4}}+\frac{2}{\sqrt{g}}\frac{\partial\sqrt{g}}{\partial w}\frac{\partial^{3}\Phi}{\partial w^{3}}+\left[\left(\frac{1}{\sqrt{g}}\frac{\partial\sqrt{g}}{\partial w}\right)^{2}+2\frac{\partial}{\partial w}\left(\frac{1}{\sqrt{g}}\frac{\partial\sqrt{g}}{\partial w}\right)\right]\frac{\partial^{2}\Phi}{\partial w^{2}}\\
\qquad\qquad+\left[\triangle_{T}\left(\frac{1}{\sqrt{g}}\frac{\partial\sqrt{g}}{\partial w}\right)+\left(\frac{\partial^{2}}{\partial w^{2}}+\frac{1}{\sqrt{g}}\frac{\partial\sqrt{g}}{\partial w}\frac{\partial}{\partial w}\right)\frac{1}{\sqrt{g}}\frac{\partial\sqrt{g}}{\partial w}\right]\frac{\partial\Phi}{\partial w}\label{eq:lap-lapl-2}
\end{multline}
where the coefficient of the first derivative $\frac{\partial\Phi}{\partial w}$
is the complete (spatial) Laplacian of $\frac{1}{\sqrt{g}}\frac{\partial\sqrt{g}}{\partial w}$. 

\bibliographystyle{plain}
\bibliography{refs}

\end{document}